\documentclass[12pt]{amsart}
\usepackage{subfigure}
\usepackage{graphicx}
\usepackage{rotating}
\usepackage{morefloats}
\usepackage[utf8]{inputenc}
\usepackage[T1]{fontenc}
\usepackage{amssymb}
\usepackage{amsmath}    
\def\R{I\kern -0,37 em R}
\def\P{I\kern -0,37 em P}
\def\Z{I\kern -0,37 em Z}

\begin{document}
\title[group invariance]{Group invariance of integrable \\ Pfaffian systems}
\author{A. Kumpera}
\address[Antonio Kumpera]{Campinas State University, Campinas, SP, Brazil}
\email{antoniokumpera@hotmail.com}

\date{August, 2016}


\maketitle

\begin{abstract}
Let $\mathcal{S}$ be an integrable Pfaffian system. If it is invariant under a transversally free infinitesimal action of a finite dimensional real Lie algebra \textit{g} and consequently invariant under the local action of a Lie group \textit{G}, we show that the \textit{vertical} variational cohomology of $\mathcal{S}$ is equal to the Lie algebra cohomology of \textit{g} with values in the space of the \textit{horizontal} cohomology in maximum dimension. This result, besides giving an effective algorithm for the computation of the variational cohomology of an invariant Pfaffian system, provides a method for detecting obstructions to the existence of finite or infinitesimal actions leaving a given system invariant.
\end{abstract}

\section{Introduction}
We study here a problem that arises naturally in connection with the integration of differential systems invariant under finite or infinitesimal group actions (\cite{Kumpera1999}) though, of course, it is also present elsewhere. The problem is to know whether a given integrable Pfaffian system $\mathcal{S}$ can be invariant under the infinitesimal action of a given Lie algebra \textit{g} or, more generally, of some Lie algebra \textit{g}. Showing the existence of such infinitesimal actions is a rather delicate problem that has to be analysed in each specific case. On the other hand, showing the non-existence can be achieved, as is common, by displaying some obstructions via cohomological methods and this is actually our main concern in this paper. Since we want to compare the variational cohomology of $\mathcal{S}$ with the Lie algebra cohomology of \textit{g}, we shall only be concerned with transversally free infinitesimal actions as defined in the section 6.

\vspace{5 mm}
\noindent
The Euler-Lagrange variational complex associated to an integrable Pfaffian system $\mathcal{S}$ is finite. As is usual, we call \textit{horizontal} that part of the complex preceding the Euler operator \textit{E} and \textit{vertical} that part subsequent to this operator. The horizontal part is a finite augmentation of \textit{E} and the vertical part a finite resolution. We show (Theorem 1) that if $\mathcal{S}$ is invariant under a transversally free infinitesimal action of the Lie algebra \textit{g} then the above finite resolution is equivalent, in positive dimensions, to the Lie algebra complex of \textit{g} taking values in the horizontal cohomology of maximum dimension. In particular, the resulting cohomology spaces are equal hence any discrepancy between the two will single out obstructions to the existence of such an infinitesimal action. The above equivalence also provides an effective method of calculating the vertical variational cohomology when the Pfaffian system is invariant.
 
\vspace{5 mm}
\noindent
Throughout the years, several authors have given distinct though essentially equivalent formulations to the variational complex. We adopt here the approach described in \cite{Kumpera1991} since it emphasizes the relationship of this complex with the algebra of generalized symmetries. Inasmuch as the usual \textit{de Rham} complex, on a manifold \textit{M}, is the differential complex associated to the algebra of all the vector fields on \textit{M}, the horizontal part of the variational complex is a de Rham complex associated to the algebra of all trivial symmetries (total derivatives) and the vertical part is also a de Rham complex associated to the algebra of all the generalized symmetries. Our first task, in this paper, consists in writing down explicitly the complex we shall be dealing with namely the restriction, to an integrable Pfaffian system, of the general complex defined in \cite{Kumpera1991}. This, unfortunately, is a rather long and boring task hence we only state, in the section 2, a well known lemma that provides all the necessary technical information relevant to the restriction procedure and thereafter construct directly, in the sections 3, 4 and 5, the desired complex. It turns out that the trivial symmetries become simply the vector fields annihilated by $\mathcal{S}$ and the generalized symmetries become the equivalence classes of the infinitesimal automorphisms of $\mathcal{S}$ modulo the trivial symmetries. There is of course nothing new about this restricted complex, just a different make-up. Invariant systems are examined in the section 6, a double complex is discussed in the section 7 and some examples are placed in evidence in the section 8. The present article is a new and updated version of \cite{Almeida2003} and the first few sections have been rewritten so as to be better suited for their application in latter sections.
 
\vspace{5 mm}
\noindent
For simplicity, we assume that all the data is $C^{\infty}$ smooth though, in each specific case,  $C^k$ smoothness for some \textit{k} will suffice. We also assume that all the manifolds are finite dimensional, connected and second countable though not necessarily orientable and that all the objects such as functions, vector fields, differential forms, etc., are globally defined on these manifolds unless stated otherwise (\textit{e.g.}, local coordinates, local generators of Pfaffian systems and distributions, etc.).

\section{Pfaffian Systems and their Prolongations}
Let $\mathcal{S}$ be a Pfaffian system, for the time being not necessarily integrable, and $\Sigma=\mathcal{S}^{\perp}$ the corresponding distribution, both defined on the manifold \textit{M} (\cite{Kumpera1999}, sect. 2). We set $dim~M=n$ and $rank~\mathcal{S}=n-p$ hence $dim~\Sigma=p,$ the last notation indicating the dimension of each contact element $\Sigma_x$ occasionally referred to as the \textit{rank} of $\Sigma.$ In terms of partial differential equations, it seems preferable to view the distribution $\Sigma$ or the Pfaffian system $\mathcal{S}$ as a section of the Grassmannian bundle $\textbf{G}^p_1M$ of $p-$dimensional linear contact elements, $p=rank~\Sigma,$ or, still better, as the sub-manifold $\mathcal{R}$ of $\textbf{G}^p_1M$ image of this section (\cite{Kumpera1991}, p. 614). Since the dimension \textit{p} will remain unchanged throughout the present discussion, we abbreviate $\textbf{G}^p_kM$ by $\textbf{G}_k,$ where \textit{k} is the order of the contact elements under consideration. Let $\mathcal{R}_1\subset\textbf{G}_2$ denote the first prolongation of $\mathcal{R}.$

\vspace{2 mm}
\newtheorem{lagrange}[LemmaCounter]{Lemma}
\begin{lagrange}
The distribution $\Sigma$ is integrable (involutive) if and only if $\mathcal{R}_1$ projects onto $\mathcal{R}.$
\end{lagrange}

\vspace{2 mm}
\noindent
A proof of this result can be found in \cite{Kumpera1999}, section 13. This proof tells us, in particular, that the fibre $(\mathcal{R}_1)_X,$ $X\in\mathcal{R},$ is either empty or else contains a single element, say \textit{Y}. The first order linear holonomic contact element $\Sigma^{(1)}_X$ associated to \textit{Y} is the unique holonomic element tangent to $\mathcal{R}$ at the point \textit{X}. Since it does not make much sense to define the variational complex for other than \textit{formally integrable} equations (or at least for equations that, after prolongation, become formally integrable at large enough orders), we see that in the present situation it becomes natural to assume that $\Sigma$ is integrable. This does not mean, however, that variational complexes cannot be associated to non-integrable Pfaffian systems. In this latter context, we ought to specify or determine the dimension $q\leq p$ of the linear contact elements for which a \textit{sufficient number} of integral contact elements do exist (\textit{e.g.}, Pfaffian systems that are \textit{in involution}, in the sense of Cartan, at the dimension \textit{q}) and consider the variational complex in the realm of the bundles $\textbf{G}^q_kM.$ We shall nevertheless restrict our attention to integrable systems.
 
\vspace{5 mm}
\noindent
When $\mathcal{S}$ (or $\Sigma$) is integrable, then $\mathcal{R}_1$ is the set of all second order contact elements determined by the $r-$dimensional integral manifolds of $\mathcal{S}$ and the assignment

\begin{equation*}
Y\in\mathcal{R}_1~\mapsto~\Sigma^{(1)}_X\subset T_X\mathcal{R},\hspace{5 mm}X=\rho_{1,2}(Y),
\end{equation*}

\vspace{3 mm}
\noindent
is an integrable distribution $\Sigma^{(1)}$ defined on the manifold $\mathcal{R}$ and equivalent to $\Sigma$ via the diffeomorphism $\beta_1=\rho_{0,1}:\mathcal{R}~\longrightarrow~M.$ In general, the $k-$th prolongation $\mathcal{R}_k\subset\textbf{G}_{k+1}$ is the set of all the $(k+1)-$st order contact elements determined by the $r-$dimensional integral manifolds of $\mathcal{S}.$ Furthermore, the assignment

\begin{equation*}
Y\in\mathcal{R}_k~\mapsto~\Sigma^{(k)}_X\subset T_X\mathcal{R}_{k-1},\hspace{5 mm}X=\rho_{k,k+1}Y,
\end{equation*}

\vspace{3 mm}
\noindent
where $\Sigma^{(k)}_X$ indicates the linear holonomic contact element at order \textit{k} associated to \textit{Y}, is an integrable distribution $\Sigma^{(k)}$ defined on the manifold $\mathcal{R}_{k-1}$ and equal to the annihilator of the restriction, to $\mathcal{R}_{k-1},$ of the canonical contact structure $\mathcal{S}_k$ of $\textbf{G}_k$ ($[\Sigma^{(k)}]^{\perp}=\iota^*\mathcal{S}_k,~\iota:\mathcal{R}_{k-1}~\hookrightarrow~\textbf{G}_k,$ \textit{cf.} \cite{Kumpera1991}, sect. 2). For any pair of integers $h\leq k,$ the distributions $\Sigma^{(k)}$ and $\Sigma^{(h)}$ are equivalent via the diffeomorphism $\rho_{h,k}:\mathcal{R}_{k-1}~\longrightarrow~\mathcal{R}_{h-1}$ ($\mathcal{R}_0=\mathcal{R}$ and $\Sigma^{(0)}=\Sigma$).

\section{The Variational Complex}
Pfaffian systems invariant under infinitesimal Lie algebra actions or, equivalently, under local Lie group actions are present in many a context and, in particular, find their place in the integration of invariant differential systems (\cite{Kumpera1999}). Motivated by some pending matters in the aforementioned paper, our concern here is to display obstructions to the existence of transversally free infinitesimal actions of a finite dimensional real Lie algebra \textit{g} leaving invariant a given integrable Pfaffian system $\mathcal{S}$ defined on a manifold \textit{M}. Exhibiting such obstructions brings some further insight into the global applicability of Lie and Cartan's methods and is achieved by comparing the variational cohomology of $\mathcal{S}$ with the Lie algebra cohomology of \textit{g} with values in an appropriate space.
 
\vspace{5 mm}
\noindent
Let $\mathcal{A}$ be the algebra of all the infinitesimal automorphisms (symmetries) of $\mathcal{S}$ and $\mathcal{B}$ the set of all the vector fields annihilated by the local sections of $\mathcal{S}$ \textit{i.e.}, the sections of $\Sigma.$ Then $\mathcal{B},$ the trivial symmetries, is an ideal of $\mathcal{A}$ and one shows, based on the Lemma 8 of \cite{Kumpera1999}, that the algebra $\mathcal{S(R)}$ of generalized symmetries of the the first order partial differential equation $\mathcal{R}$ associated to $\mathcal{S}$ ($S(D)$ in the notations of \cite{Kumpera1991}, sect. 9) identifies with $\mathcal{A/B}.$ Further, $\mathcal{A}$ is a module over the ring $\mathcal{I}$ of all the first integrals of $\mathcal{S}$ and $\mathcal{B}$ is a module over the ring $\mathcal{F}$ of all the functions on \textit{M}, hence $\mathcal{S(R)}$ is a module over $\mathcal{I}.$ We adopt here the approach to the variational complex as described in \cite{Kumpera1991} for it emphasizes its relationship with the algebra $\mathcal{S(R)}$ and seems best suited for our purposes. In what follows, we describe what the general (full) horizontal and vertical complexes as well as the Euler operator, considered in \cite{Kumpera1991}, become upon restriction to the equation $\mathcal{R}.$   

\section{The Horizontal Operator}
We now construct directly the \textit{horizontal} part of the variational complex associated to an integrable Pfaffian system $\mathcal{S}$ namely, that part preceding the Euler operator. Let us denote, as above, by $\mathcal{B}$ the algebra of all the (globally defined) vector fields $\eta$ tangent to the distribution $\Sigma=\mathcal{S}^{\perp}$ \textit{i.e.}, $\eta\in\Gamma(\Sigma),$ and by $\mathcal{F}$ the ring of all the (globally defined) $C^{\infty}$ functions on the underlying manifold \textit{M}. The dual space $\mathcal{H}=\mathcal{B}^*,$ with respect to the $\mathcal{F}-$module structure, is equal to the set of all the global sections of the dual bundle $\Sigma^*\simeq T^*M/\mathcal{S}$ and, correspondingly, $\wedge\mathcal{H}\simeq\Gamma(\wedge\Sigma^*).$ The differential

\begin{equation*}
d_H:\wedge^s\mathcal{H}~\longrightarrow~\wedge^{s+1}\mathcal{H}
\end{equation*}

\vspace{3 mm}
\noindent
is defined by the usual formula:

\begin{equation}
d_H\mu(\eta_1,~\cdots~,\eta_{s+1})=\sum_i~(-1)^{i+1}\theta(\eta_i)\mu(\eta_1,~\cdots~,\widehat{\eta_i},~\cdots~,\eta_{s+1})~+
\end{equation}
\begin{equation*}
\hspace{23 mm}\sum_{i<j}~(-1)^{i+j}\mu([\eta_i,\eta_j],\eta_1,~\cdots~,\widehat{\eta_i},~\cdots~,\widehat{\eta_j},~\cdots~,\eta_{s+1}),
\end{equation*}

\vspace{3 mm}
\noindent 
where the $\eta_i$ are vector fields tangent to $\Sigma,~\theta(\eta_i)$ is the usual Lie derivative and $[\eta_i,\eta_j]$ is the usual Lie bracket. Let $(U;x^1,~\cdots~,x^p,y^1,~\cdots~,y^q).~q=n-p.$ be a foliated chart of $\Sigma$ for which the integral manifolds, in \textit{U}, are given by the equations $y^{\lambda}=c^{\lambda},$ the $c^{\lambda}$ being constants. Then, an element of $\wedge^s\mathcal{H}$ has the local expression

\begin{equation*}
\mu=\sum~a_{i_1\cdots i_s}~dx^{i_1}\wedge~\cdots~\wedge dx^{i_s},
\end{equation*}

\vspace{3 mm}
\noindent
the coefficients $a_{i_1\cdots i_s}$ being $C^{\infty}$ functions on \textit{U}, and

\begin{equation}
d_H\mu=\sum~(da_{i_1\cdots i_s}|\Sigma)~dx^{i_1}\wedge~\cdots~\wedge dx^{i_s}
\end{equation}
\begin{equation*}
\hspace{14 mm}=\sum~\frac{\partial a_{i_1\cdots i_s}}{\partial x^i}~dx^i\wedge dx^{i_1}\wedge~\cdots~\wedge dx^{i_s},
\end{equation*}

\vspace{3 mm}
\noindent
where $da_{i_1\cdots i_s}|\Sigma$ (resp. $dx^i$) stands for the restriction of this differential to the integral manifolds of $\Sigma$ or, still more simply, its restriction to the contact elements $\Sigma_x,~x\in M.$

\vspace{5 mm}
\noindent
Let us now extend the differential (1) by adding, in the cochains, a term that corresponds in \cite{Kumpera1991} to the module $\mathcal{C}$ of all the contact $1-$forms. Here we consider $\mathcal{C}=\Gamma(\mathcal{S})$ to be the module of all the global sections of $\mathcal{S},$ take the cochain space $\Phi^{r,s}=(\wedge^r\mathcal{C})\otimes_{\mathcal{F}}(\wedge^s\mathcal{H})$ and consider its elements as horizontal forms with values in $\wedge^r\mathcal{C}$. The extended differential 

\begin{equation*}
d_H:\wedge^r\mathcal{C}\otimes\wedge^s\mathcal{H}~\longrightarrow~\wedge^r\mathcal{C}\otimes\wedge^{s+1}\mathcal{H}
\end{equation*}

\vspace{3 mm}
\noindent
is then defined by

\begin{equation}
d_H(\omega\otimes\mu)(\eta_1,~\cdots~,\eta_{s+1})=\sum_i~(-1)^{i+1}\theta(\eta_i)[\mu(\eta_1,~\cdots~,\widehat{\eta_i},~\cdots~,\eta_{s+1})\omega]~+
\end{equation}
\begin{equation*}
\hspace{13 mm}\sum_{i<j}~(-1)^{i+j}\mu([\eta_i,\eta_j],\eta_1,~\cdots~,\widehat{\eta_i},~\cdots~,\widehat{\eta_j},~\cdots~,\eta_{s+1})\omega,
\end{equation*}

\vspace{3 mm}
\noindent
where $\eta_i\in\Gamma(\Sigma)$ and $\theta(\eta_i)$ is the Lie derivative. The second term, on the right hand side, belongs of course to $\wedge^r\mathcal{C},$ the same being true for the first term since $\mathcal{S}$ is integrable and consequently $\theta(\eta_i)\mathcal{C}\subset\mathcal{C}.$ Let $(U;x^1,~\cdots~,x^p,y^1,~\cdots~,y^q)$ be a foliated chart for the distribution $\Sigma$ in which the integral manifolds are given by the slices $y^{\lambda}=c^{\lambda}.$ Then a typical element of $\Phi^{r,s}$ is locally a sum of terms

\begin{equation*}
\mu=a~dy^{j_1}\wedge~\cdots~\wedge dy^{j_r}\otimes dx^{i_1}\wedge~\cdots~\wedge dx^{i_s},
\end{equation*}

\vspace{3 mm}
\noindent
and

\begin{equation}
d_H\mu=\sum~\frac{\partial a}{\partial x^i}~dy^{j_1}\wedge~\cdots~\wedge dy^{j_r}\otimes dx^i\wedge dx^{i_1}\wedge~\cdots~\wedge dx^{i_s},
\end{equation}

\vspace{3 mm}
\noindent
where $dx^i$ stands for the restriction $dx^i|\Sigma.$ The last formula as well as the formula (2), though helpful in theoretical considerations, is most often useless in practice since it requires the local integration of $\Sigma.$ We can nevertheless remedy this situation as follows. Let us consider any coordinate system $(U;x^i,y^j)$ with the sole requirement that the family $\{dx^i|\Sigma\}$ be free at every point of \textit{U}, thus providing a field of co-frames for $\Sigma^*|U.$ Next, we consider the local basis $\{\eta_i\}$ of $\Sigma$ defined by the conditions $<\eta_i,dx^j>~=\delta^j_i.$ Since $\Sigma$ is integrable and since each $\eta_i$ projects onto $\partial/\partial x^i,$ it follows that $[\eta_i,\eta_j]=0.$ A typical element of $\Phi^{r,s}$ can now be written locally as a sum of terms

\begin{equation*}
\mu=a~\omega^{j_1}\wedge~\cdots~\wedge \omega^{j_r}\otimes dx^{i_1}\wedge~\cdots~\wedge dx^{i_s},
\end{equation*}

\vspace{3 mm}
\noindent
where $dx^i$ stands for $dx^i|\Sigma$ and $\{\omega^{\lambda}=dy^{\lambda}-\sum~Y^{\lambda}_i dx^i\}$ is a local basis of $\mathcal{S}.$ The formula (3) then reduces to 

\begin{equation}
d_H\mu=\sum_i~[\theta(\eta_i)(a~\omega^{j_1}\wedge~\cdots~\wedge \omega^{j_r})]\otimes dx^i\wedge dx^{i_1}\wedge~\cdots~\wedge dx^{i_s},
\end{equation}

\vspace{3 mm}
\noindent
and a similar formula can also replace (2), the derivatives $\partial a_{\cdots}/\partial x^i$ being replaced by $\theta(\eta_i)(a_{\cdots}).$

\vspace{5 mm}
\noindent
Let $\mathcal{I}^r=\mathcal{I}^r(\mathcal{S})$ denote the module of all the globally defined \textit{invariant forms}
$\omega$ of degree \textit{r} with respect to the Pfaffian system $\mathcal{S}$ namely, those forms satisfying the condition (\cite{Kumpera1999}, sect. 4):

\begin{equation*}
\theta(\eta)\omega=0,\hspace{5 mm}\forall~\eta\in\Gamma(\mathcal{S}^{\perp}).
\end{equation*}

\vspace{3 mm}
\noindent
Then, $\mathcal{I}^0=\mathcal{I}$ is the ring of all the global first integrals of $\mathcal{S},~\oplus_r~\mathcal{I}^r$ is a graded $\mathcal{I}-$subalgebra of $\mathcal{A}$ and the formula (3) shows that the \textit{short} sequence 

\begin{equation}
0~\longrightarrow~\mathcal{I}^r~\longrightarrow~\Phi^{r,0}~\xrightarrow{d_H}~\Phi^{r,1}
\end{equation}

\vspace{3 mm}
\noindent
is exact for each integer \textit{r}. We next show that the \textit{long} sequence

\begin{equation}
0~\longrightarrow~\mathcal{I}~\longrightarrow~\Phi^{0,0}~\xrightarrow{d_H}~\Phi^{0,1}~\xrightarrow{d_H}~\cdots~\xrightarrow{d_H}~\Phi^{0,p}~\longrightarrow~0
\end{equation}

\vspace{3 mm}
\noindent
is locally exact. In fact, let $(U;x^i,y^{\lambda})$ be a foliated chart for $\Sigma.$ Then the formula (2) defines, for each fixed set of values $y^{\lambda}=c^{\lambda},$ the differential of the \textit{de Rham} complex on the corresponding slice, whereupon results the local exactness of (7) since the usual homotopy operators can be written incorporating the parameters $y^{\lambda}.$ Let us finally show that

\begin{equation}
\Phi^{r,p-1}~\xrightarrow{d_H}~\Phi^{r,p}~\longrightarrow~0
\end{equation}

\vspace{3 mm}
\noindent
is locally exact. A typical element of $\Phi^{r,p}$ is, locally, a sum

\begin{equation*}
\omega=\sum~a_{j_1\cdots j_r}~dy^{j_1}\wedge~\cdots~\wedge dy^{j_r}\otimes dx^{i_1}\wedge~\cdots~\wedge dx^{i_p},
\end{equation*}

\vspace{3 mm}
\noindent
hence, upon integrating for example along $x^1,$ we obtain the element

\begin{equation*}
\Omega=\sum~A_{j_1\cdots j_r}~dy^{j_1}\wedge~\cdots~\wedge dy^{j_r}\otimes dx^{i_2}\wedge~\cdots~\wedge dx^{i_p},\hspace{2 mm}\frac{\partial A_{[j]}}{\partial x^1}=a_{[j]}
\end{equation*}

\vspace{3 mm}
\noindent
such that $d_H\Omega=\omega.$

\section{The Vertical Operator}
Let us next construct the \textit{vertical} part of the Euler-Lagrange complex namely, that part subsequent to the Euler operator. We recall that $\mathcal{A}$ denotes the algebra of all the infinitesimal automorphisms of $\mathcal{S}$ and $\mathcal{B}$ the ideal of those vector fields that are tangent to $\Sigma=\mathcal{S}^{\perp}.$ The system $\mathcal{S}$ being integrable, any vector field tangent to $\Sigma$ is an infinitesimal automorphism. Based on the Lemma 1, it can be shown as observed earlier, that the algebra $\mathcal{S(R)}$ of the generalized symmetries of the equation $\mathcal{R}$ associated to $\mathcal{S}$ identifies with $\mathcal{A/B}.$ The following remarks will be used later: 

\vspace{5 mm}
(a) $\mathcal{A}$ is a module over the ring $\mathcal{I}$ and $\mathcal{B}$ is a module over $\mathcal{F}$ hence $\mathcal{S(R)}$ is a module over $\mathcal{I}.$ 

\vspace{3 mm}
(b) When $\xi\in\mathcal{A}$ is tangent to an integral leaf of $\mathcal{S}$ at some point $x_0,$ then it is also tangent to that leaf at all the other points.

\vspace{5 mm}
\noindent
We now define $\Xi^r=\Phi^{r,s}/d_H\Phi^{r,s-1},$ denote by $\mathbf{q}_r:\Phi^{r,s}~\longrightarrow~\Xi^r$ the quotient map and observe that $\Xi^r$ is an $\mathcal{I}-$module since $d_H\mathcal{I}=0.$ An element $\omega\otimes\mu\in(\wedge^r\mathcal{C})\otimes(\wedge^s\mathcal{H})$ can also be considered as an $\mathcal{I}-$multilinear form on $\mathcal{S(R)}$ with values in $\wedge^s\mathcal{H}$ by setting

\begin{equation*}
(\omega\otimes\mu)([\xi_1],~\cdots~,[\xi_r])=\omega(\xi_1,~\cdots~,\xi_r)\mu,\hspace{5 mm}\xi_i\in\mathcal{A},
\end{equation*}

\vspace{3 mm}
\noindent
where each $[\xi_i]$ denotes the class of $\xi_i~mod~\mathcal{B}.$ In fact, when $\eta\in\mathcal{B},$ then $(\omega\otimes\mu)(~\cdots~,\eta,~\cdots~)=0$ and consequently $(\omega\otimes\mu)([\xi_1],~\cdots~,[\xi_r])$ is well defined on $\mathcal{S(R)}.$ Furthermore (\cite{Kumpera1991}, sect. 11), since 

\begin{equation*}
[d_H(\omega\otimes\mu)](\xi_1,~\cdots~,\xi_r)=d_H[(\omega\otimes\mu)(\xi_1,~\cdots~,\xi_r)],
\end{equation*}

\vspace{3 mm}
\noindent
the form $d_H(\omega\otimes\mu),~\omega\otimes\mu\in(\wedge^r\mathcal{C})\otimes(\wedge^s\mathcal{H}),$ considered as a multilinear form on $\mathcal{S(R)},$ takes values that vanish under the projection $\mathbf{q}_0:\Phi^{0,s}~\longrightarrow~\Xi^0.$ Hence, to any element $\tau\in\Xi^r,$ we can associate an $\mathcal{I}-$multilinear form $[\tau],$ defined on $\mathcal{S(R)}$ and taking values in $\Xi^0,$ as follows. We take $\Omega=\sum~\omega_i\otimes\mu_i\in\Phi^{r,s}$ such that $\mathbf{q}_r(\Omega)=\tau$ and set  

\begin{equation*}
[\tau]([\xi_1],~\cdots~,[\xi_r])=\mathbf{q}_0\Omega([\xi_1],~\cdots~,[\xi_r]).
\end{equation*}

\vspace{3 mm}
\noindent
The mapping $\tau\mapsto[\tau]$ being injective (\cite{Kumpera1991}, sect.11), we are led to consider the defining formula

\begin{equation*}
d_V(\mathbf{q}_r(\omega\otimes\mu))([\xi_1],~\cdots~,[\xi_{r+1}])=
\end{equation*}
\begin{equation}
=\mathbf{q}_0\{\sum_i~(-1)^{i+1}\theta([\xi_i])[\omega([\xi_1],~\cdots~,\widehat{[\xi_i]},~\cdots~,[\xi_{r+1}])\mu]\}+
\end{equation}
\begin{equation*}
\mathbf{q}_0\{\sum_{i<j}~(-1)^{i+j}\omega([[\xi_i],[\xi_j]],[\xi_1],~\cdots~,\widehat{[\xi_i]},~\cdots~,\widehat{[\xi_j]},~\cdots~,[\xi_{r+1}])\mu\},
\end{equation*}

\vspace{3 mm}
\noindent
where $\omega\otimes\mu\in\Phi^{r,s}$ and $\xi_i\in\mathcal{A}.$ The second term on the right hand side clearly belongs to $\Xi^0.$ As for the first term, let us write 

\begin{equation*}
\mu_i=\omega([\xi_1],~\cdots~,\widehat{[\xi_i]},~\cdots~,[\xi_{r+1}])\mu\in\Phi^{0,s}
\end{equation*}

\vspace{3 mm}
\noindent
and let us assume that some $\xi_i=\eta\in\mathcal{B}.$ If $j\neq i,$ then $\mu_i=0$ and if $j=i$ then, since $d_H\mu_i=0,$

\begin{equation*}
\theta(\eta)\mu_i=i(\eta)d_H\mu_i+d_Hi(\eta)\mu_i=d_Hi(\eta)\mu_i
\end{equation*}

\vspace{3 mm}
\noindent
and consequently $\mathbf{q}_0\theta(\eta)\mu_i=0.$ In any case,, (9) provides a well defined multilinear form on $\mathcal{S(R)}$ taking values in $\Xi^0$ and it can be shown (\textit{e.g.}, in coordinates) that the multilinear form $d_V(\mathbf{q}_r(\omega\otimes\mu))$ is the image of an element $\tau\in\Xi^{r+1}.$

\vspace{5 mm}
\noindent
The Euler-Lagrange variational complex associated to the integrable Pfaffian system $\mathcal{S}$ is the finite sequence ($p=rank~\Sigma,~q=rank~\mathcal{S}$)

\begin{equation}
0~\longrightarrow~\mathcal{I}~\longrightarrow~\Phi^{0,0}~\xrightarrow{~d_H~}~\Phi^{0,1}~\xrightarrow{~d_H~}~\cdots~\xrightarrow{~d_H~}~\Phi^{0,p-1}~\xrightarrow{~d_H~}~
\end{equation}
\begin{equation*}
\hspace{9 mm}\xrightarrow{~d_H~}~\Phi^{0,p}~\xrightarrow{~E~}~\Xi^1~\xrightarrow{~d_V~}~\cdots~\xrightarrow{~d_V~}~\Xi^{q-1}~\xrightarrow{~d_V~}~\Xi^q~\longrightarrow~0,
\end{equation*}

\vspace{3 mm}
\noindent
where \textit{E}, the \textit{Euler} operator, is the composite 

\begin{equation*}
\Phi^{0,p}~\xrightarrow{~\mathbf{q}_0~}~\Xi^0~\xrightarrow{~d_V~}~\Xi^1
\end{equation*}

\vspace{3 mm}
\noindent
and $q=rank~\mathcal{S}.$ This complex is locally exact and reduces, locally, to (7) since $\Xi^r$ vanishes on account of (8). 

\section{Invariant Pfaffian Systems}
Let $\mathcal{S}$ be an integrable Pfaffian system invariant under the infinitesimal action $\Phi:g~\longrightarrow~\chi(M)$ of a finite dimensional real Lie algebra \textit{g} (\cite{Kumpera1999}, sect. 3). For every $v\in g,~\Phi(v)\in\mathcal{A}$ hence the action induces a Lie algebra morphism $\Phi:g~\longrightarrow~\mathcal{S(R)}.$ If, further, we assume that the action $\Phi$ is transversally free (\cite{Kumpera1999}, sect. 4) \textit{i.e.}, if

\vspace{5 mm}
$~(i)~dim~g=dim~\Phi(g)_x,~\forall~x\in M,~and$

\vspace{3 mm}
$(ii)~T_xM=\Sigma_x\oplus\Phi(g)_x,~\forall~x\in M,$

\vspace{5 mm}
\noindent
then the above morphism is injective and the following result holds:

\vspace{2 mm}
\newtheorem{malgrange}[LemmaCounter]{Lemma}
\begin{malgrange}
The algebra $\mathcal{S(R)}$ is generated by $\Phi(g)$ over the ring $\mathcal{I}$ of first integrals of $\mathcal{S},$ any $\mathbf{R}-$basis of $\Phi(g)$ is an $\mathcal{I}-$basis of $\mathcal{S(R)}$ and the $\mathcal{I}-$dual $\mathcal{S(R)}^*$ identifies with $\mathcal{I}^1.$
\end{malgrange}

\vspace{2 mm}
\noindent
Let us next recall some properties of the invariant forms associated to $\mathcal{S}.$ By definition (\cite{Kumpera1999}, sext. 4), the exterior form $\omega$ is an invariant form of $\mathcal{S}$ if $\theta(\eta)\omega=0$ for all the vector fields $\eta\in\Gamma(\Sigma),~\Sigma=\mathcal{S}^{\perp}.$ It follows that $\theta(f\eta)\omega=0$ for any function \textit{f}, hence $\omega$ is invariant if and only if $i(\eta)\omega=i(\eta)d\omega=0,$ for all $\eta\in\Gamma(\Sigma).$ We infer that a necessary and sufficient condition for $\omega$ to be an invariant form is that it be expressible, locally, in terms of the first integrals of $\mathcal{S}$ and their differentials. When $\omega$ is invariant then so are the forms $f\omega,~d\omega$ and $\theta(\xi)\omega,$ where \textit{f} is a first integral and $\xi$ an infinitesimal automorphism of the system $\mathcal{S},$ hence the set of all the invariant forms is a differential algebra over the ring $\mathcal{I}$ invariant under the infinitesimal action $\Phi$ via the Lie derivative. Let $\{v_i\}$ be a basis of \textit{g}. Then the linear forms $\omega^i\in\mathcal{C}$ defined by the conditions $<\Phi(v_i),\omega^j>=\delta^j_i$ are a global basis of invariant forms of $\mathcal{S},$ a so-called Cartan basis (\cite{Kumpera1999}, sect. 6,8), and

\begin{equation}
d\omega^i=\sum_{j<k}~c^i_{jk}\omega^j\wedge\omega^k,
\end{equation}

\vspace{3 mm}
\noindent
where $\{-c^i_{jk}\}$ is the set of structure constants of \textit{g} with respect to the above chosen basis. The real sub-space $\Omega\subset\Gamma(T^*M)$ generated by the forms $\omega^i$ only depends on $\Phi$ and acts as an $\mathbf{R}-$dual to the space $h=\Phi(g).$ Let us denote by $\mathcal{F}$ the ring of $C^{\infty}$ functions on \textit{M} and by $\wedge\Omega$ the exterior algebra of $\Omega$ over the field $\mathbf{R}.$ Since the set $\{\omega^i\}$ is a global basis of the Pfaffian system $\mathcal{S},$ it follows that $\mathcal{C}\simeq\Omega\otimes_{\mathbf{R}}\mathcal{F}$ and, more generally, that

\begin{equation}
\Phi^{r,s}=(\wedge^r\mathcal{C})\otimes_{\mathcal{F}}(\wedge^s\mathcal{H})\simeq(\wedge^r\Omega)\otimes_{\mathbf{R}}(\wedge^s\mathcal{H})\simeq(\wedge^r\mathcal{S(R)}^*)\otimes_{\mathcal{I}}(\wedge^s\mathcal{H}).
\end{equation}

\vspace{3 mm}
\noindent
Furthermore, since $d\omega(\xi_1,\xi_2)=-\omega([\xi_1,\xi_2]),$ for any $\omega\in\Omega$ and $\xi_i\in\Phi(g),$ it also follows that the formula (9) reduces, whenever $\omega\in\wedge^r\Omega$ and $\xi_i\in\Phi(g)\subset\mathcal{S(R)},$ to the expression

\begin{equation}
[d\omega\otimes\mathbf{q}_0\mu+(-1)^{deg~\omega}~\omega\wedge d_V(\mathbf{q}_0\mu)](\xi_1,~\cdots~,\xi_{r+1}),
\end{equation}

\vspace{3 mm}
\noindent
where $\mathbf{q}_0\mu(\xi)$ is the Lie derivative $\mathbf{q}_0(\theta(\xi)\mu).$ We shall see later that the above formula still holds for any invariant form $\omega$ since $d_H\omega=0$ implies $d_V\omega=d\omega.$

\vspace{5 mm}
\noindent
We next consider the elements $\tau\in\Xi^r$ as $\mathcal{I}-$multilinear (and skew-symmetric) forms $[\tau]$ defined on $\mathcal{S(R)}$ and taking values in $\Xi^0$ (\textit{cf.} sect. 4). Each form $[\tau]$ restricts to an $\mathcal{R}-$multilinear form $\varrho=[\tau]_{\mathbf{R}}$ defined on the real sub-space $h\subset\mathcal{S(R)}$ and conversely, on account of the Lemma 2, each $\mathcal{R}-$multilinear form $\varrho$ defined on \textit{h} and taking values in $\tau\in\Xi^0$ extends, by $\mathcal{I}-$multilinearity, to a form $[\tau]$ defined on $\mathcal{S(R)}$ and verifying $[\tau]_{\mathbf{R}}=\varrho.$ Furthermore, since in the realm of real vector spaces the sub-space $d_H(\Phi^{0,s-1})$ admits a complement in $\Phi^{0,s},$ the form $\varrho$ lifts to a real form $\overline{\varrho}$ defined on \textit{h} and taking values in $\Phi^{0,s}=\wedge^s\mathcal{H}.$ Finally, making use of $\mathcal{F}-$multilinearity on \textit{h}, $\overline{\varrho}$ extends to an $\mathcal{F}-$multilinear form $\tilde{\varrho}$ defined on $\chi(M)$ and with values in $\wedge^s\mathcal{H}$ by requiring that $i(\eta)\tilde{\varrho}=0$ for all the vector fields $\eta\in\Gamma(\Sigma).$ We thus obtain an element $\tilde{\varrho}\in\Phi^{r,s}$ such that $[\tilde{\varrho}]_{\mathbf{R}}=\varrho$ and, consequently, the assignment $\tau~\longmapsto~[\tau]$ of the section 4 becomes bijective when $\mathcal{S}$ is invariant under a transversally free infinitesimal action. Identifying \textit{g} with \textit{h} and thereafter $\Omega$ with $g^*,$ the expression (13) or, equivalently, the formula (9) shows that the cochain complex defining the cohomology of \textit{g} with values in $\Xi^0$ and relative to the representation $\rho(v)(\mathbf{q}_0\mu)=\mathbf{q}_0(\theta(\xi)\mu),$ with $v\in g,~\mu\in\Phi^{0,s}$ and $\xi=\Phi(v),$ is equal to

\begin{equation}
\Xi^0~\xrightarrow{~d_V~}~\Xi^1~\xrightarrow{~d_V~}~\Xi^2~\xrightarrow{~d_V~}~\cdots~.
\end{equation}

\vspace{3 mm}
\noindent
Since $\mathbf{q}_0$ is surjective, we can rewrite the above complex by

\begin{equation}
\Phi^{0,p}~\xrightarrow{~E~}~\Xi^1~\xrightarrow{~d_V~}~\Xi^2~\xrightarrow{~d_V~}~\cdots~,
\end{equation}

\vspace{3 mm}
\noindent
without affecting the cohomology groups in positive dimensions, the latter being the vertical part of the Euler-Lagrange complex associated to $\mathcal{S}$ namely, the finite resolution of \textit{E}.

\vspace{2 mm}
\newtheorem{orange}[TheoremCounter]{Theorem}
\begin{orange}
Let $\mathcal{S}$ be an integrable Pfaffian system invariant under a transversally free infinitesimal action of the Lie algebra $g.$ Then the finite resolution of the Euler operator $E$ is equal, in positive dimensions, to the cochain complex of the Lie algebra $g$ with values in $\Xi^0.$
\end{orange}

\vspace{2 mm}
\noindent
Since $E=d_V\circ\mathbf{q}_0,$ we also infer that the space of cocycles in $\Phi^{0,p}$ (\textit{i.e.}, $ker~E$) is equal to the inverse image, by $\mathbf{q}_0,$ of the $0-$dimensional cohomology of \textit{g} namely, the inverse image of the sub-space composed by the $g-$invariant elements of $\Xi^0$ ($\rho(v)(\mathbf{q}_0\mu)=0,~\forall v\in g$). Given an integrable Pfaffian system $\mathcal{S},$ we can now confront the variational cohomology with the cohomology of \textit{g} taking values in $\Xi^0$ and eventually detect obstructions to the existence of a transversally free infinitesimal action of \textit{g} leaving the system $\mathcal{S}$ invariant. When $dim_{\mathbf{R}}~\Xi^0<\infty,$ which is very seldom the case, then we usually have more information on the cohomology of \textit{g}. For instance, when \textit{g} is semi-simple, its cohomology vanishes in dimensions one and two (Whitehead's lemmas) and consequently the same must hold for the variational cohomology.

\vspace{5 mm}
\noindent
The above discussion leads us finally to define a \textit{twisted complex} very convenient in calculations. Identifying the algebra \textit{g} with its image $h=\Phi(g),$ the elements of $g^*$ can be considered as those belonging to a sub-set of the invariant forms considered previously and we define the complex $(\wedge g^*\otimes_{\mathbf{R}}\Xi^0,\partial)$ by setting:

\begin{equation}
\partial(\omega\otimes\mathbf{q}_0\mu)=d\omega\otimes\mathbf{q}_0\mu+(-1)^{deg~\omega}~\omega\wedge d_V\mathbf{q}_0\mu,
\end{equation}

\vspace{3 mm}
\noindent
$d\omega$ being the differential of the invariant form $\omega$ and $(d_V\mathbf{q}_0\mu)(\xi)$ the Lie derivative $\mathbf{q}_0(\theta)(\xi)\mu.$ The inclusion $g\hookrightarrow\mathcal{S(R)}$ determines, by transposition, a natural transformation (cochain complex morphism) from the Euler-Lagrange complex towards the twisted chain complex associated to \textit{g} and this, in turn, induces morphisms on the corresponding cohomology groups. If, further, the structure constants of \textit{g} are known, then the first term on he right hand side of (16) is easily calculable. We also consider the trivial extension $(\wedge g^*\otimes_{\mathbf{R}}\Xi^0,d\otimes Id)$ of the standard cochain complex associated to \textit{g}. The corresponding cohomology is clearly isomorphic to the cohomology of \textit{g} and it becomes very useful to compare, via homological methods, the twisted complex with the standard one (\textit{e.g.}, homotopical deformation of $d\otimes Id$ into $\partial$). 

\section{A double complex}
Usually, the Euler-Lagrange complex is defined in terms of a double complex. We shall also define, in this section, a double complex not for the sake of re-defining what has already been defined earlier but in view of obtaining further invariants for the Pfaffian system $\mathcal{S}.$ We begin by stating an almost obvious result that will however be needed in the sequel. Two first order distributions (field of contact elements) $\Sigma_V$ and $\Sigma_H,$ defined on a manifold \textit{M}, will be said to be \textit{complementary} when $T_xM=(\Sigma_V)_x\oplus(\Sigma_H)_x$ for any $x\in M,$ which means that the tangent bundle $TM$ is the direct sum of the two distributions. Let us denote by \textit{V} and \textit{H} the corresponding projections on each factor and by $d_V$ and $d_H$ the associated \textit{Frölicher-Nijenhuis} differentials of type \textit{d} (\cite{Frolich1956}). Recalling some definitions, given a vector valued $1-$form $u:TM\longrightarrow TM$ defined on \textit{M} (\textit{e.g.}, the projections \textit{V} and \textit{H}) and a scalar differential form $\omega$ of degree \textit{r}, we shall denote by $i_u\omega$ the scalar form (indicated in the previous reference by $u\overline{\wedge}\omega$) defined as follows:

\begin{equation}
i_u\omega(\xi_1,~\cdots~,\xi_r)=\sum_j~\omega(\xi_1,~\cdots~,\xi_{j-1},u(\xi_j),\xi_{j+1},~\cdots~,\xi_r).
\end{equation}

\vspace{3 mm}
\noindent
The map $i_u:\wedge M^*~\longrightarrow~\wedge M^*$ becomes an interior derivation of degree zero namely, it verifies the identity

\begin{equation}
i_u(\omega\wedge\mu)=(i_u\omega)\wedge\mu+\omega\wedge i_u\mu
\end{equation}

\vspace{3 mm}
\noindent
and \textit{interior} means that $i_u$ vanishes on functions. The associated Frölicher-Nijenhuis derivation $d_u=[i_u,d]=i_u\circ d-d\circ i_u$ is of degree one and

\begin{equation}
d_u(\omega\wedge\mu)=(d_u\omega)\wedge\mu+(-1)^{deg~\omega}~\omega\wedge d_u\mu.
\end{equation}

\vspace{3 mm}
\noindent
Moreover, the map $u\mapsto d_u$ is $\mathbf{R}-$linear, $d_{Id}=d$ and $[d_u,d]=d_u\circ d+d\circ d_u=0$ \textit{i.e.}, $d_u$ is a derivation of type \textit{d}. In particular, $d_H+d_V=d$ since $H+V=Id.$ We now cite the desired result.

\vspace{2 mm}
\newtheorem{solange}[TheoremCounter]{Theorem}
\begin{solange}
The distributions $\Sigma_V$ and $\Sigma_H$ are simultaneously integrable if and only if $d_V^2=0$ or, equivalently, if and only if $d_H^2=0,$ where $\Sigma_V$ and $\Sigma_H$ are arbitrary complementary distributions and V, H the corresponding projections. 
\end{solange}

\vspace{2 mm}
\noindent
Let us now assume that both distributions are integrable. Then

\begin{equation}
0=[d_{V+H},d_{V+H}]=[d_V,d_V]+[d_H,d_H]+[d_V,d_H]+[d_H,d_V]
\end{equation}

\vspace{3 mm}
\noindent
and we infer that $0=[d_V,d_H]+[d_H,d_V]=2(d_Vd_H+d_Hd_V),$ where after $d_Vd_H=-d_Hd_V.$ Let $\Psi^{r,s}$ denote the module of all the differential forms of type $(r,s)$ with respect to the splitting $TM=\Sigma_V\oplus\Sigma_H.$ We shall say that a vector field $\xi$ is vertical (resp. horizontal) whenever it is tangent to $\Sigma_V$ (resp. $\Sigma_H$) and observe that a differential form $\omega$ of degree $r+s$ belongs to $\Psi^{r,s}$ if and only if $\omega(\xi_1,~\cdots~\xi_{r+s})=0$ whenever more than \textit{r} entries are vertical or more than \textit{s} entries are horizontal. Forms of type $(r,0)$ are those forms that vanish whenever some entry is horizontal and, similarly, for the forms of type $(0,s),$ this property translating formally by $i(\xi)\omega=0,~\xi\in\Gamma(\Sigma_H),$ and correspondingly for the type $(0,s)$ forms. Since all the modules of differential forms here considered are projective of finite type, we infer that $\Psi^{r,s}=\Psi^{r,0}\wedge\Psi^{0,s},$ $\Psi^{r,0}=\wedge^r\Psi^{1,0},$ $\Psi^{0,s}=\wedge^s\Psi^{0,1},$ $\Psi^{1,0}=\mathcal{C}$ and $\Psi^{0,1}$ is the module of all the forms vanishing on $\Sigma_V.$ We now show that $\Psi^{0,1}\simeq\mathcal{H}.$ In fact, any $\omega\in\Psi^{0,1}$ induces, by restriction, an element of $\mathcal{H}$ and, conversely, any element of $\mathcal{H}$ extends (lifts) to an element of $\Psi^{0,1}$ by requiring that it vanishes on the distribution $\Sigma_V.$ More generally, the elements of $\Psi^{0,s}$ are in one-to-one correspondence with the elements of $\wedge^s\mathcal{H},$ the mapping  $\Psi^{0,s}\longrightarrow\wedge^s\mathcal{H}$ being simply the restriction of the forms to $\Sigma_H,$ with the inverse defined by $\mu\mapsto\mu\circ\wedge^rH$ (read \textit{the projection H}). We can also identify $\Psi^{1,0}=\mathcal{C}$ and, more generally, $\Psi^{r,0}=\wedge^r\mathcal{C}$ with the modules $\Gamma(\Sigma^*_V)$ and $\Gamma(\wedge^r\Sigma^*_V)$ respectively. Given any two complementary integrable distributions, they both play technically similar roles and can be, for whatever purpose, interchanged. Nevertheless, our earlier approach to the Euler-Lagrange complex, privilege was accorded to the initially given distribution $\Sigma.$ We now identify $\Phi^{r,s}$ with $\Psi^{r,s}$  by firstly identifying $\wedge^s\mathcal{H}$ with $\Phi^{0,s}$ and secondly by identifying the elements $\omega\otimes\mu\in\Phi^{r,s}$ with $\omega\wedge\mu\in\Psi^{r,s}.$ Under this identification, the operator $d_H,$ defined in (3), identifies with the Frölicher-Nijenhuis operator $(-1)^rd_H$ and, furthermore, we can now proceed to write down the desired double complex ($p=dim~\Sigma_H,~q=dim~\Sigma_V$ and $d_H$ is equal to (3) or to $(-1)^rd_H$ in the Frölicher-Nijenhuis sense). Obviously, we shall choose $\Sigma_H=\Sigma$ but the first component $\Sigma_V$ has no privileged choices, the resulting double complexes being equivalent for any two different choices. Nevertheless, in order to obtain the largest benefit with respect to the additional invariants, we shall choose $\Sigma_V=\Phi(g)~i.e.$, $(\Sigma_V)_x=(\Phi(g))_x,$ for all $x\in M,$ since the Lie algebra action provides considerable additional information. Let us write: 

\vspace{5 mm}

\begin{equation*}
\hspace{9 mm}0\hspace{18 mm}0\hspace{32 mm}0\hspace{20 mm}0    
\end{equation*}
\begin{equation*}
\hspace{9 mm}\downarrow\hspace{18 mm}\downarrow\hspace{32 mm}\downarrow\hspace{20 mm}\downarrow
\end{equation*}
\begin{equation*}
\hspace{9 mm}\mathcal{J}^0(\mathcal{S})\xrightarrow{~d~}\mathcal{J}^1(\mathcal{S})\xrightarrow{~d~}\cdots\xrightarrow{~d~}\mathcal{J}^{p-1}(\mathcal{S})\xrightarrow{~d~}\mathcal{J}^p(\mathcal{S})  
\end{equation*}
\begin{equation*}
\hspace{10 mm}\downarrow\hspace{18 mm}\downarrow\hspace{32 mm}\downarrow\hspace{20 mm}\downarrow
\end{equation*}
\begin{equation*}
0\rightarrow\mathcal{I}^0(\mathcal{S})\rightarrow\Phi^{0,0}\hspace{2 mm} \xrightarrow{d_H}\hspace{2 mm}\Phi^{0,1}~\xrightarrow{d_H}~\cdots~\xrightarrow{d_H}~\Phi^{0,p-1}~\xrightarrow{d_H}~\Phi^{0,p}~\xrightarrow{q_0}~\Xi^0
\end{equation*}
\begin{equation*}
\hspace{10 mm}d\downarrow\hspace{7 mm}d_V\downarrow\hspace{11 mm}d_V\downarrow\hspace{25 mm}d_V\downarrow\hspace{13 mm}d_V\downarrow\hspace{6 mm}d_V\downarrow\hspace{3 mm}
\end{equation*}
\begin{equation*}
0\rightarrow\mathcal{I}^1(\mathcal{S})\rightarrow\Phi^{1,0}\hspace{2 mm}\xrightarrow{d_H}\hspace{2 mm}\Phi^{1,1}~\xrightarrow{d_H}~\cdots~\xrightarrow{d_H}~\Phi^{1,p-1}~\xrightarrow{d_H}~\Phi^{1,p}~\xrightarrow{q_1}~\Xi^1   
\end{equation*}
\begin{equation*}
\hspace{10 mm}d\downarrow\hspace{7 mm}d_V\downarrow\hspace{11 mm}d_V\downarrow\hspace{25 mm}d_V\downarrow\hspace{13 mm}d_V\downarrow\hspace{6 mm}d_V\downarrow\hspace{3 mm}   
\end{equation*}
\begin{equation*}
\hspace{14 mm}\vdots\hspace{15 mm}\vdots\hspace{19 mm}\vdots\hspace{33 mm}\vdots\hspace{21 mm}\vdots\hspace{14 mm}\vdots\hspace{6 mm}    
\end{equation*}
\begin{equation*}
\hspace{14 mm}\downarrow\hspace{14 mm}\downarrow\hspace{18 mm}\downarrow\hspace{32 mm}\downarrow\hspace{20 mm}\downarrow\hspace{13 mm}\downarrow\hspace{10 mm}    
\end{equation*}
\begin{equation*}
0\rightarrow\mathcal{I}^r(\mathcal{S})~\rightarrow\Phi^{r,0}\hspace{2 mm}\xrightarrow{d_H}\hspace{2 mm}\Phi^{r,1} ~\xrightarrow{d_H}~\cdots~\xrightarrow{d_H}~\Phi^{r,p-1}\hspace{2 mm}\xrightarrow{d_H}~\Phi^{r,p}~\xrightarrow{q_r}~\Xi^r
\end{equation*}
\begin{equation*}
\hspace{11 mm}d\downarrow\hspace{7 mm}d_V\downarrow\hspace{11 mm}d_V\downarrow\hspace{25 mm}d_V\downarrow\hspace{13 mm}d_V\downarrow\hspace{6 mm}d_V\downarrow\hspace{3 mm}    
\end{equation*}
\begin{equation*}
~0\rightarrow\mathcal{I}^q(\mathcal{S})~\rightarrow\Phi^{q,0}\hspace{2mm}\xrightarrow{d_H}\hspace{2mm}\Phi^{q,1}~\xrightarrow{d_H}~\cdots~\xrightarrow{d_H}~\Phi^{q,p-1}\hspace{2 mm}\xrightarrow{d_H}~\Phi^{q,p}~\xrightarrow{q_q}~\Xi^q     
\end{equation*}

\vspace{5 mm}
\noindent
where $q_{\ell}=\mathbf{q}_{\ell},~\forall \ell,$ in the previous notation and, on the one before the last line, $r$ stands for $q-1$ (due to the lack of space). The operators $d_V:\Xi^{\ell}\longrightarrow\Xi^{\ell+1}$ are those defined by the formula (9) and the relation $d_H\circ d_V=-d_V\circ d_H,$ in the Frölicher-Nijenhuis sense, now yields 

\begin{equation*}
(-1)^{r+1}d_H\circ d_V=(-1)^rd_V\circ d_H,
\end{equation*}

\vspace{3 mm}
\noindent
hence all the operators in the above diagram commute. We already know that the lines of this diagram are exact at $\Phi^{\ell,0},$ $\mathcal{I}^{\ell}(\mathcal{S})$ being the module of the invariant forms of degree $\ell$ with respect to the system $\mathcal{S}.$ Let us now compute the kernel of $\Phi^{0,\ell}\longrightarrow\Phi^{1,\ell}.$ We first observe that a form $\omega\in\Psi^{r,s}$ is equal to zero if and only if $\omega(\xi_1,~\cdots~,\xi_r,\eta_1,~\cdots~,\eta_s)=0$ for all the entries \textit{r} of which are vertical and \textit{s} are horizontal. Furthermore, since vertical vector fields are spanned by $\Phi(g),$ it will suffice to consider only the elements $\xi_i\in\Phi(g).$ With this in mind and recalling that $[\xi,\eta]\in\Gamma(\Sigma_H)$ whenever $\xi\in\Phi(g)$ and $\eta\in\Gamma(\Sigma_H),$ a straightforward calculation leads to the formula

\begin{equation*}
d_V\mu(\xi,\eta_1,~\cdots~,\eta_{\ell})=[\theta(\xi)\mu](\eta_1,~\cdots~,\eta_{\ell}),
\end{equation*}

\vspace{3 mm}
\noindent
where $\mu\in\Psi^{0,\ell},$ $\xi\in\Phi(g)$ and $\eta_i\in\Gamma(\Sigma_H).$ Consequently, $d_V\mu=0$ if and only if $\theta(\xi)\mu=0$ for all $\xi\in\Phi(g).$ The space $\mathcal{J}^{\ell}(\mathcal{S}),$ kernel of $d_V:\Phi^{0,\ell}\longrightarrow\Phi^{1,\ell},$ is equal to the $\mathcal{I}(\mathcal{S})-$module of the horizontal forms of degree $\ell$ invariant under the infinitesimal action $\Phi.$ Since the relation $[\xi,\eta]\in\Gamma(\Sigma_H)$ continues to hold for $\eta$ horizontal and $\xi$ an arbitrary infinitesimal automorphism of $\Sigma_H$ tangent to $\Sigma_V,$ all the elements of $\mathcal{J}^{\ell}(\mathcal{S})$ are invariant under such automorphisms. Still better, the relation 

\begin{equation*}
\theta(f\xi)\mu=f\theta(\xi)\mu+df\wedge i(\xi)\mu
\end{equation*}

\vspace{3 mm}
\noindent
shows that all the elements of $\mathcal{J}^{\ell}(\mathcal{S})$ are invariant by any vertical vector field (tangent to $\Sigma_V$). It should however be remarked that the Lie derivative $\theta(\xi)\mu$ only has a sense, along infinitesimal automorphisms $\xi$ of $\Sigma_H,$ when $\Phi^{0,\ell}$ is considered in its original definition as $\wedge^{\ell}\mathcal{H}.$ A similar calculation will also show, once more, that the kernel of $d_H:\Phi^{\ell,0}\longrightarrow\Phi^{\ell,1}$ is the set of all the invariant forms of degree $\ell$ with respect to $\Sigma_H.$ Since $d=d_V+(-1)^rd_H$ ($d_H$ in the sense of the previous diagram) and $d_H(\mathcal{I}^{\ell}(\mathcal{S}))=0,$ we obtain the complex

\begin{equation}
0\longrightarrow\mathbf{R}\longrightarrow\mathcal{I}^0(\mathcal{S})\xrightarrow{~d~}\mathcal{I}^1(\mathcal{S})\xrightarrow{~d~}~\cdots~\xrightarrow{~d~}\mathcal{I}^q(\mathcal{S})\longrightarrow 0
\end{equation}

\vspace{3 mm}
\noindent
whose cohomology is the invariant cohomology of $\Sigma~(=\Sigma_H).$ It only depends on $\Sigma$ since the cochains are the invariant forms. When $\Sigma$ is simple (\textit{i.e.}, when \textit{M} admits a differentiable quotient structure, modulo the integral leaves of $\Sigma$) then (21) is equivalent to the \textit{de Rham} complex of the quotient manifold $M/\Sigma.$
We can thus read on (21) the non-quotientability of $\Sigma$ or, at least, discard prospective candidates to quotients by scanning the cohomological invariants (\textit{e.g.}, Betti and Lefschetz numbers). Since in our context of a transversally free $g-$action we are provided with a global basis of invariant $1-$forms (\textit{cf.} the section 6), the cohomology of (21) can, in principle, be calculated as soon as we know the constants of structure of \textit{g}. Similarly, we can also display a complex

\begin{equation}
0\longrightarrow\mathbf{R}\longrightarrow\mathcal{J}^0(\mathcal{S})\xrightarrow{~d~}\mathcal{J}^1(\mathcal{S})\xrightarrow{~d~}~\cdots~\xrightarrow{~d~}\mathcal{J}^p(\mathcal{S})\longrightarrow 0
\end{equation}

\vspace{3 mm}
\noindent
whose cochains are the $g-$invariant horizontal forms. Its cohomology is the equivariant cohomology of $(\Sigma,\Phi).$ Our earlier remark shows as well that (22) is also the complex of $g-$invariant forms with respect to $\Sigma_V.$ When the integral foliation of the distribution $\Sigma_V$ admits a differentiable quotient, this complex is equivalent to the \textit{de Rham} complex of the quotient manifold $M/\Sigma_V.$ It should further be emphasized that neither (21) nor (22) do rely upon the construction of the double complex exhibited on the previous diagram. 

\vspace{5 mm}
\noindent
There are other invariants of $\Sigma$ that can be detected on the double complex. For instance, we can detect the \textit{relatively invariant} forms of $\Sigma$ ($\omega$ is said to be relatively invariant when its differential $d\omega$ is invariant). Assuming that $\omega$ is of degree $\ell,$ let us write

\begin{equation*}
\omega=\omega^{\ell,0}+\omega^{\ell-1,1}+\omega^{\ell-2,2}+~\cdots~+\omega^{1,\ell-1}+\omega^{0,\ell},
\end{equation*}

\vspace{3 mm}
\noindent
where $\omega^{r,s},~r+s=\ell$ denotes the component of type $(r,s)$ of $\omega$ and, using $d_H$ in the sense of the double complex, we can also write

\begin{equation*}
d\omega=d_V\omega^{\ell,0}+(-1)^{\ell}d_H\omega^{\ell,0}+d_V\omega^{\ell-1,1}+(-1)^{\ell-1}d_H\omega^{\ell-1,1}+d_V\omega^{\ell-2,2} +
\end{equation*}
\begin{equation*}
 +d_V\omega^{\ell-2,2}+(-1)^{\ell-2}d-H\omega^{\ell-2,2}+~\cdots~+d_V\omega^{1,\ell-1}-d_H\omega^{1,\ell-1}+d_V\omega^{0,\ell}+d_H\omega^{0,\ell}.   
\end{equation*}

\vspace{3 mm}
\noindent
If $d\omega$ is invariant then it belongs to $\mathcal{I}^{\ell+1}(\mathcal{S})\subset\Phi^{\ell+1,0}.$ Writing out, symbolically, the types in the expression of $d\omega,$ we find the following display:

\begin{equation*}
d\omega=(\ell+1,0)+(\ell,1)+(\ell,1)+(\ell-1,2)+(\ell-1,2)+(\ell-2,3)+~\cdots~+
\end{equation*}
\begin{equation*}
 +(2,\ell-1)+(1,\ell)+(1,\ell)+(0,\ell+1),   
\end{equation*}

\vspace{3 mm}
\noindent
hence, in order that $d\omega$ be of type $(\ell+1,0),$ it is necessary and sufficient that all the remaining terms cancel out namely, the following equalities must be verified:

\begin{equation}
(-1)^{\ell}d_H\omega^{\ell,0}+d_V\omega^{\ell-1,1}=0\hspace{10 mm}(-1)^{\ell-1}d_H\omega^{\ell-1,1}+d_V\omega^{\ell-2,2}=0
\end{equation}
\begin{equation*}
\cdots\hspace{10 mm}-d_H\omega^{1,\ell-1}+d_V\omega^{0,\ell}=0\hspace{10 mm}d_H\omega^{0,\ell}=0.
\end{equation*}

\vspace{3 mm}
\noindent
Let us now assume that the above conditions do hold. Then $d\omega$ is of type $(\ell+1,0)$ and

\begin{equation*}
d_Hd\omega=d_Hd_V\omega^{\ell,0}=-d_Vd_H\omega^{\ell,0}=d_V((-1)^{\ell}d_V\omega^{\ell-1,1})=0
\end{equation*}

\vspace{3 mm}
\noindent
on account of the first relation in (23), hence $d\omega$ is invariant. The conditions (23) can be retraced in the spectral sequence associated to the double complex and with respect to the \textit{vertical} filtration. At present, we shall not inquire any further in these matters since they escape our main concern namely, the Euler-Lagrange complex that appears at the upper right corner of the double complex, the composite $d_V\circ q_0=q_1\circ d_V$ being the Euler operator. In our notations, $\mathcal{J}^p(\mathcal{S})$ is the set of null Lagrangians. In particular, when $\omega=\omega^{\ell,0},$ then (23) reduces to $d_H\omega^{\ell,0}=0$ and $d\omega$ is invariant if and only if $\omega$ is invariant.

\vspace{5 mm}
\noindent
Let us now return to the end of the section 6. We consider a Pfaffian system $\mathcal{S}$ defined on the manifold \textit{M} and invariant under the transversally free infinitesimal action $\Phi:g\longrightarrow\chi(M),$ as well as a second system $\widetilde{\mathcal{S}},$ defined on $M\times N,$ constructed via a basis of invariant forms of $\mathcal{S}$ and a second  infinitesimal action $\Psi:g\longrightarrow\chi(N).$ In order to avoid any confusion, we denote by $\overline{\mathcal{S}}$ the system, on $M\times N,$ shifted of $\mathcal{S}$ namely, the system  generated by $p^*_M\mathcal{S}+p^*_NT^*N,$ where $p_M$ and $p_N$ are the projections onto the respective factors. This system is invariant under the \textit{horizontal} action $\overline{\Phi}:g\longrightarrow\chi(M\times N)$ defined by $\overline{\Phi}(v)=(\Phi(v),0).$ We also consider the \textit{vertical} action $\widetilde{\Phi}:g\longrightarrow\chi(M\times N)$ given by $\widetilde{\Phi}(v)=(0,\Psi(v)).$ Finally, we set $\overline{\Sigma}=\overline{\mathcal{S}}^{\perp}$ and $\widetilde{\Sigma}=\widetilde{\mathcal{S}}^{\perp}.$ Then $\overline{\Sigma}\subset\widetilde{\Sigma},$ $\widetilde{\mathcal{S}}\subset\overline{\mathcal{S}}$ and the inclusion $\iota:\overline{\Sigma}\hookrightarrow\widetilde{\Sigma}$ transposes to $\iota^*:\Gamma(\wedge\widetilde{\Sigma}^*)\longrightarrow\Gamma(\wedge\overline{\Sigma}^*).$ Since vector fields tangent to $\overline{\Sigma}$ are also tangent to $\widetilde{\Sigma},$ we infer, from the formula (1), that $d_H\circ\iota^*=\iota^*\circ d_H,$ hence $\iota^*$ is a natural transformation of cochain complexes. The same, however, is not true for the variational parts of the Euler-Lagrange complexes associated to $\overline{\mathcal{S}}$ and $\widetilde{\mathcal{S}}$ since $rank~\overline{\Sigma}\neq rank~\widetilde{\Sigma}$ (unless $\Sigma=TM$). Define $\overline{\Phi}\oplus\widetilde{\Phi}:g\oplus g\longrightarrow\chi(M\times N)$ by setting, in the obvious way, $\overline{\Phi}\oplus\widetilde{\Phi}(v,w)=\overline{\Phi}(v)+\widetilde{\Phi}(w).$ In studying this product infinitesimal action that is transversally free with respect to the sum of the aforementioned distributions, we shall be able to find new invariants for the system $\mathcal{S}.$ Nevertheless, we shall not pursue any further this line of investigation since, for the time being, we do not foresee any interesting applications for these new invariants.

\section{Examples}
Throughout this section, we replace integrable Pfaffian systems by the corresponding integral foliations since the latter are much more visible. Though all the foliations considered in the sequel are entirely naive, the resulting homological calculations are not always so.

\vspace{5 mm}
\noindent
\textbf{Example 1}\hspace{3 mm}\textit{The torus.}

\vspace{5 mm}
\noindent
Let $\mathcal{F}$ be the foliation on the $2-$dimensional torus \textit{T} whose leaves are the cosets of a $1-$dimensional sub-group \textit{H}. Then $\mathcal{F}$ is invariant under the infinitesimal action generated by any element of the Lie algebra of \textit{T} ($\equiv \mathbf{R}^2$), this action being transversally free when the chosen element does not belong to the Lie algebra \textit{h} of \textit{H}. The variational cohomology at $\Xi^1$ is equal to $\mathbf{R}$ and, using Green's formula, it can be shown that the cohomology class of an element $[\omega]\in\Xi^1$ identifies with the real number $\int_T\omega.$ When the slope of the element in \textit{h} is rational, the calculations are very simple and both spaces $\Xi^0$ and $\Xi^1$ identify with the set of all the global first integrals of $\mathcal{F}.$ However, when this slope is irrational, the global first integrals of $\mathcal{F}$ reduce to the constants (the curves become everywhere dense) and it becomes more involved to describe the spaces $\Xi^0$ and $\Xi^1$ as well as to calculate directly the cohomology. This is an example where the advantages of the Lie algebra cohomology calculations become apparent. 

\vspace{5 mm}
\noindent
\textbf{Example 2}\hspace{3 mm}\textit{The Möbius strip.}

\vspace{5 mm}
\noindent
Let $\mathcal{F}$ be the foliation on the Möbius strip $\mathcal{M}$ whose leaves are the "double" circles, except for the central circle (under the usual identification $(1,y)\equiv(-1,-y),$ where $\mathcal{F}$ is the foliation induced by the segments parallel to the $x-$axis). Both spaces $\Xi^0$ and $\Xi^1$ identify again with the set of all the global first integrals of $\mathcal{F}$ that, in turn, identifies with the set of all the even functions defined on the interval $]-1,1[.$ The variational cohomology at $\Xi^1$ vanishes and, whatever the representation $\rho:g=\mathbf{R}~\longrightarrow~Der~\Xi^0,$ the cohomology of \textit{g} in dimension one is non-trivial (the derivative of an even function usually ceases to be even). The above disagreement shows that $\mathcal{F}$ cannot be invariant by any transversally free infinitesimal action, a fact that is geometrically obvious since such an action would provide an orientation to $\mathcal{M}.$

\vspace{5 mm}
\noindent
\textbf{Example 3}\hspace{3 mm}\textit{Foliation with a compact attractor.}

\vspace{5 mm}
\noindent
We consider, on the infinite cylinder $\mathcal{C}=\mathbf{R}\times\mathcal{S}^1$ with the coordinates $(t,\theta),$ the foliation $\mathcal{F}$ obtained by integrating the vector field

\begin{equation}
\eta=t\frac{\partial}{\partial t}+\frac{\partial}{\partial\theta}~.
\end{equation}

\vspace{3 mm}
\noindent
The nature of the spaces $\Xi^0$ and $\Xi^1$ is rather involved but a straightforward calculation shows that the cohomology at $\Xi^1$ is null. On the other hand, whatever the representations of the Lie algebra $g=\mathbf{R}$ into $Der~\Xi^0,$ the Lie algebra cohomology, in dimension one, cannot vanish. Consequently, the foliation $\mathcal{F}$ does not admit any $1-$dimensional transversally free infinitesimal action that leaves it invariant. This fact is also geometrically obvious since the local $1-$parameter group $(\phi_u)$ generated by any such infinitesimal action would be defined, for small \textit{u}, on a whole neighborhood of the limit circle $\{0\}\times\mathcal{S}^1$ and would transform this circle into open compact subsets of the neighboring leaves, this being of course excluded.

\vspace{5 mm}
\noindent
An entirely similar situation arises in the double solid torus (two solid tori glued together by their boundaries) upon taking the Reeb foliation inside each of the tori. The common boundary torus is the unique compact leaf.

\vspace{5 mm}
\noindent
\textbf{Example 4}\hspace{3 mm}\textit{Spheres and rays.}

\vspace{5 mm}
\noindent
On the space $M=\mathbf{R}^{p+1}-\{0\},$ let $\mathcal{F}_1$ be the foliation whose leaves are the spheres centered at the origin and $\mathcal{F}_2$ the foliation whose leaves are the rays issued from the origin. We first consider the \textit{sphere} foliation $\mathcal{F}_1$ and calculate $\Xi^0,$ one possible argument being as follows: Each element $\mu\in\Phi^{0,q}$ identifies canonically with a differentiable $1-$parameter family $(\overline{\mu}_{\rho})$ of differentiable $q-$forms defined on the unit sphere $\mathcal{S}^p$ and, under this identification, $d_H\mu\in\Phi^{0,q+1}$ also identifies with $d\overline{\mu}_{\rho}.$ We next take a differentiable $1-$parameter family $(\overline{\mu}_t),~t>0,$ of $p-$forms on $\mathcal{S}^p.$ Then, upon choosing a fixed volume form $\Omega$ on this unit sphere (\textit{e.g.}, the volume form associated to the induced Euclidean metric), we can determine, upon integration, a differentiable function $\varphi:\mathbf{R}_+~\longrightarrow~\mathbf{R}$ such that $\overline{\mu}_t-\varphi(t)\Omega$ is, for each \textit{t}, a coboundary. Restating the Lemma 4.2 of \cite{Sternberg1964} in its stronger version (as is proved in the subsequent two pages), we can use it to establish a stronger $1-$parameter version of the mentioned Lemma 4.2 and prove in the aforementioned context that there exists a differentiable $1-$parameter family $\overline{\eta}_t$ of $(p-1)-$forms defined on $\mathcal{S}^p$ such that $\overline{\mu}_t-\varphi(t)\Omega=d\overline{\eta}_t.$ Returning to $\Phi^{0,p}$ and taking the form $\widetilde{\Omega}=r^*\Omega,~r:X~\mapsto~(1/\parallel X\parallel)X,$ defined on \textit{M}, we conclude that each $\mu\in\Phi^{0,p}$ determines a differentiable function $\varphi$ such that $\mu-\varphi\widetilde{\Omega}=d_H\eta,$ where $\eta\in\Phi^{0,p-1},$ and consequently that $\Xi^0$ is equal to the set of all the real-valued differentiable functions defined on $\mathbf{R}_+$ \textit{i.e.}, equal to the set of all the global first integrals of $\mathcal{F}_1.$

\vspace{5 mm}
\noindent
Let us now calculate $\Xi^1$. Observing that $d\rho$ is a global generator of the Pfaffian system that annihilates $\mathcal{F}_1,$ any element $\eta_1\in\Phi^{1,p}$ writes $\eta_1=\eta\wedge d\rho,$ with $\eta\in\Phi^{0,q},$ and $d_H\eta_1=(d_H\eta)\wedge d\rho$ since $d_H(d\rho)=0.$ Consequently, the element $\mu_1=\mu\wedge d\rho\in\Phi^{1,p}$ is equal to $d_H\eta_1,$ with $\eta_1\in\Phi^{1,p-1},$ if and only if $\mu=d_H\eta$ hence the present calculation reduces to the previous one and $\Xi^1$ is again equal to the set of all the global first integrals of $\mathcal{F}_1.$ 

\vspace{5 mm}
\noindent
It now becomes easy to show that the variational cohomology at $\Xi^1$ is null. The foliation $\mathcal{F}_1$ is of course invariant under many $1-$dimensional transversally free infinitesimal actions and the vanishing of the Lie algebra cohomology in dimension one can easily be verified.

\vspace{5 mm}
\noindent
We next consider the \textit{radial} foliation $\mathcal{F}_2.$  Here we can procede locally, on open sets saturated by rays, and integration along these rays will show that $\Xi^q=0$ for $0\leq q\leq p.$ The variational as well as the Lie algebra cohomologies vanish, their comparison not revealing the following geometrical facts:

\vspace{5 mm}
(\textit{a}) When \textit{p} is even, there cannot exist a transversally free infinitesimal action leaving the radial foliation $\mathcal{F}_2$ invariant. In fact, since the tangent spaces to $\mathcal{F}_1$ and $\mathcal{F}_2$ are complementary, any such infinitesimal action would project onto the spheres producing an infinitesimal action operating tangentially to the spheres and, in restriction to these spheres, would be free. However, even dimensional spheres do not admit nowhere vanishing vector fields.

\vspace{5 mm}
(\textit{b}) When \textit{p} is odd, such transversally free infinitesimal actions do exist only for $p=1,3.$ Their non-existence for $p=7$ is essentially a consequence of the fact that $\mathcal{S}^7$ is not a group manifold and, for all the other values of \textit{p}, that the corresponding spheres are not parallelizable.

\vspace{5 mm}
\noindent
We can enhance the variational cohomology by adding non-trivial cocycles to the space \textit{M}. For example, let us take for the manifold \textit{M} the portion of $\mathbf{R}-\{0\}$ in between the spheres $\mathcal{S}^p(1)$ and $\mathcal{S}^p(2),$ (the numbers indicating the radii) and identify these two spheres by the radial map. Then $\mathcal{F}_1$ induces a foliation $\overline{\mathcal{F}}_1$ in spheres, $\mathcal{F}_2$ a foliation $\overline{\mathcal{F}}_2$ in circles (in fact, $M={S}^p\times{S}^1$) and we can show, for the foliation $\overline{\mathcal{F}}_2,$ that $\Xi^0$ is equal to the set of all the differentiable functions defined on the sphere $\mathcal{S}^p(1)$ or, equivalently, to the set of all the global first integrals of $\overline{\mathcal{F}}_2.$ Furthermore, $\Xi^r$ is equal to the product of $(^p_r)$ copies of $\Xi^0$ and  $\Xi^p=\Xi^0.$ As for the variational cohomology, we can apply again the $1-$parameter version of the Lemma 4.2 mentioned earlier and conclude that it vanishes at $\Xi^{r+1}$ whenever $r+1<p$ and that it is equal to $\mathbf{R}$ at $\Xi^p.$ Stokes' formula will then show that the cohomology class of an element $[\omega]\in\Xi^p$ identifies with the real number $\int_M\omega.$

\vspace{5 mm}
\noindent
Returning to the geometrical facts described earlier, we can retrace (\textit{a}) by looking at the variational cohomology. In fact, since any global vector field on an even dimensional sphere has a singularity, whatever the representation $\rho$ of a Lie algebra \textit{g} into $Der~\Xi^0\simeq\chi(\mathcal{S}^p),$ the corresponding Lie algebra cohomology cannot vanish in dimension one. As for the property (\textit{b}), it requires a deeper analysis that seems out of reach in the present setting. Nevertheless, it can be shown that transversally free abelian infinitesimal actions leaving $\overline{\mathcal{F}}_2$ invariant cannot exist since the corresponding Lie algebra cohomologies with values in $\Xi^0$ would vanish in dimension \textit{p} thus contradicting the variational cohomology.

\vspace{5 mm}
\noindent
\textbf{Example 5}\hspace{3 mm}\textit{Exterior differential systems.}

\vspace{5 mm}
\noindent
Exterior differential systems are a generalization of Pfaffian systems and allow the presence of exterior differential forms of arbitrary degree. Almost all that has been discussed in the present paper can be extended to these systems as soon as integrability is replaced by \textit{involutiveness}, within the appropriate definitions. These systems were introduced by Élie Cartan in the beginning of the $20^{th}$ century and afterwards much exploited by him in view of solving fundamental geometrical problems. (\cite{Cartan1945},\cite{Cartan1953},\cite{Kumpera1962}, \textit{see also} \cite{Lie1895}). We shall not discuss here this extended context since much has to be re-evaluated including the initial notion of generalized symmetry. Moreover, it might be considerably beneficial to replace Lie group actions by Lie groupoid actions since the latter will embrace many more situations.

\vspace{5 mm}
\noindent
$\mathbf{Acknowledgement.}$ Acknowledgement is due to the \textit{Klingon Foundation} for the support received.

\vspace{5 mm}
\noindent

\begin{figure}[h!]
\centering
\includegraphics[scale=3.4]{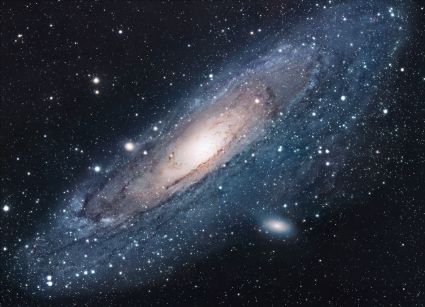}
\caption{The Universe}
\label{fig:univerise}
\end{figure}

\bibliographystyle{plain}
\bibliography{references}

\begin{thebibliography}{1}

\bibitem{Almeida2003}
R.~Almeida, A.~Kumpera, and J.~Rubin.
\newblock {On the Variational Cohomology of $g-$invariant Foliations}.
\newblock {\em J. Math. Phys.}, 44:4702--4712, 2003.

\bibitem{Cartan1945}
E.~Cartan.
\newblock {\em Les systèmes différentiels extérieurs et leurs applications
  géométri-ques}.
\newblock Hermann, Paris, 1945.

\bibitem{Cartan1953}
E.~Cartan.
\newblock {\em O\!euvres Complètes, Part. II}.
\newblock Gauthier-Villars, Paris, 1953.

\bibitem{Frolich1956}
A.~Frölicher and A.~Nijenhuis.
\newblock {Theory of vector-valued differential forms, Part I}.
\newblock {\em Indag. Math.}, 18:338--359, 1956.

\bibitem{Kumpera1962}
A.~Kumpera.
\newblock {\em {Exterior Differential Systems}}.
\newblock Princeton Notes, 1962.

\bibitem{Kumpera1991}
A.~Kumpera.
\newblock {Les symétries généralisées et le complexe d'Euler-Lagrange}.
\newblock {\em J. Fac. Sci. Univ. Tokyo}, 38:589--622, 1991.

\bibitem{Kumpera1999}
A.~Kumpera.
\newblock {On the Lie and Cartan theory of invariant differential systems}.
\newblock {\em J. Math. Sci. Univ. Tokyo}, 6:229--314, 1999.

\bibitem{Lie1895}
S.~Lie.
\newblock {Verwertung des Gruppenbegriffes für Differentialgleichungen, I}.
\newblock {\em Ber. Ges. Leipzig, Math-Phys.}, 15:261--322, 1895.

\bibitem{Sternberg1964}
S.~Sternberg.
\newblock {\em Lectures on differential geometry}.
\newblock Prentice-Hall, Inc., 1964.

\end{thebibliography}
\end{document}